\documentclass[11pt]{article}

\usepackage{latexsym}
\usepackage{epsfig}
\usepackage{float}
\usepackage{amsfonts}
\restylefloat{figure}
\restylefloat{table}
\usepackage{array}

\newcommand{\Section}[1]{\setcounter{equation}{0} \section{#1}}

\newcommand{\rf}[1]{(\ref{#1})}
\newcommand{\rfE}[1]{\emph{(\ref{#1})}}
\newcommand{\lb}{\label}

\newcommand{\beq}[1]{ \begin{equation}\label{#1} }
\newcommand{\eeq}{\end{equation}}
\newcommand{\ba}{\begin{array}}
\newcommand{\ea}{\end{array}}
\newcommand{\beqa}{\begin{eqnarray}}
\newcommand{\eeqa}{\end{eqnarray}}
\newcommand{\beqan}{\begin{eqnarray*}}
\newcommand{\eeqan}{\end{eqnarray*}}

\DeclareRobustCommand{\qed}{%
  \ifmmode 
  \else \leavevmode\unskip\penalty9999 \hbox{}\nobreak
  \fi
  \quad\hbox{\qedsymbol}}
\newcommand{\qedsymbol}{\rule{2pt}{4pt}}
{\par\addvspace{6pt}\normalfont \itshape #1\@addpunct{.}\hskip\labelsep\ignorespaces\normalfont}{\qed\par\addvspace{6pt}}

\newtheorem{thm}{Theorem}
\newtheorem{corollary}[thm]{Corollary}

\newtheorem{rem}[thm]{Remark}
\newtheorem{rems}[thm]{Remarks}
\newtheorem{prop}[thm]{Proposition}

\newcommand{\f}[1]{ \ensuremath{ {\mathbb  #1} } }

\newcommand{\sv}[2]{ {\scriptstyle s}^{#1}_{#2} }

\newcommand{\bt}[1]{~~~\mbox{#1}~~~}
\newcommand{\hs}{\hspace{.1in}}
\newcommand{\vs}{\vskip .10in}

\newcommand{\half}[1]{\frac{#1}{2}}
\newcommand{\inv}[1]{\frac{1}{#1}}

\newcommand{\e}{\emph{e}}

\newcommand{\ep}{\hfill \rule{2.4mm}{2.4mm}}

\newcommand{\is}[1]{#1^{\ast}}

\newcommand{\lm}{\lambda}

\newcommand{\vp}{\varphi}
\newcommand{\w}{\omega}

\newcommand{\q}{QMF}
\newcommand{\qs}{QMFs}

\newlength{\captsize} \let\captsize=\footnotesize
\newlength{\captwidth} 
\DeclareFontFamily{U}{msb}{}
\DeclareFontShape{U}{msb}{m}{n}{ <5> <6> <7> <8> <9> gen * msbm
               <10> <10.95> <12> <14.4> <17.28> <20.74> <24.88> msbm10}{} 
\DeclareSymbolFont{AMSb}{U}{msb}{m}{n}
\DeclareMathSymbol{\N}{\mathalpha}{AMSb}{"4E}
\DeclareMathSymbol{\R}{\mathalpha}{AMSb}{"52}
\DeclareMathSymbol{\Z}{\mathalpha}{AMSb}{"5A}

\begin{document}
%
%
%
%
%
\title{Linear Phase Perfect Reconstruction Filters
and\\ Wavelets with Even Symmetry.}
\author{Lucas  Monz\'{o}n\\
Department of Applied Mathematics\\
University of Colorado\\
Boulder, CO 80309-0526}
\date{December 1, 1999}
\maketitle
%

\begin{abstract}
Perfect reconstruction filter banks can be used to generate
a variety of wavelet bases.  Using IIR
linear phase filters one can  obtain symmetry properties for the wavelet
and scaling functions. 

In this paper we describe all possible IIR linear
phase filters generating symmetric wavelets with any prescribed number
of vanishing moments. In analogy with the well known FIR case, we
construct and study a new family of wavelets obtained by considering
maximal number of vanishing moments for each fixed order of the IIR
filter.  Explicit expressions for the coefficients of numerator,
denominator, zeroes, and poles are presented.

This new parameterization  allows one to design linear phase
quadrature
mirror filters with
many other properties of interest such as  filters that have
any preassigned set of zeroes in the stopband
or that satisfy an almost interpolating property.

Using Beylkin's approach, 
it is indicated how to implement these IIR
filters not as recursive filters but as FIR filters. 
\end{abstract}
\vs
\vs
\vs


\Section{Introduction}
\lb{int}

Quadrature mirror filters (\qs)\footnote{Originally, 
filters leading to perfect reconstruction were named conjugate
quadrature filters~\cite{SMI-BAR:1986} but we use the term \q\ as
in~\cite[Pages 162-163]{DAUBEC:1992}.}  allow one to design two-channel
perfect reconstruction filter banks and therefore to generate
wavelet bases. For implementation reasons we
either consider finite impulse response (FIR) or realizable infinite
impulse response (IIR) filters. 
The first class corresponds to polynomial
filters, the second, wider class, to rational filters.  
In both cases, we consider filters
with real coefficients and we identify
filters with their z-transforms. 
Rational \qs\ admit
properties difficult or impossible to achieve with polynomial \qs. For
example, neither linear phase (except for the Haar filter) nor
interpolation properties can be exactly obtained even
though both properties can be approximately achieved using FIR
coiflets \cite{MO-BE-HE:1999}.  In contrast, both properties can be
exactly obtained with IIR filters but not simultaneously.
Parameterizations of polynomials \cite{POLLEN:1990,VAI-HOA:1988} or rational
\qs\ \cite{DOG-VAI:1990} are well known. In both cases, linear
conditions on the coefficients of the filter, such as vanishing
moments for the wavelet or the scaling function, become
non-linear conditions on the parameters describing the coefficients.
This problem can sometimes be overcome, and there are explicit
families of rational filters with an arbitrary number of zeroes at $-1,$
that is an arbitrary number of vanishing moments for the wavelet. But
the additional requirement of symmetry conditions, which correspond to
linear phase filters, again forces us to solve non-linear
systems. For this reason, there are only some
examples of linear phase filters generating wavelets with vanishing
moments ~\cite{HER-VET:1993, SELESN:1999}.  

Often, it is much simpler to describe the filters in terms of their
absolute values (as in \cite{CAG-AKA:1993, DAUBEC:1988} for
polynomials or~\cite{HER-VET:1993} for rational \qs), leading to
implicit constructions. In most cases the solutions are obtained
numerically and the exact coefficients are not known except for small
filter lengths.  These numerical solutions are \qs\ within a certain
accuracy. In practical applications this is not an inconvenience
because we always work within a certain precision that can in
principle be attained except for large filter lengths.  Nevertheless,
we may like to incorporate the accuracy sought as a parameter in the
process of design-implementation of the \q\@.

Leaving aside the problem of accuracy once the size of the filter
increases, we come across a much deeper problem.  Namely, 
how does one
simultaneously impose different properties in an exact or approximate
sense?  A good example is Daubechies' parameterization of the
magnitude response of polynomial \qs. Even though she used it to
obtain the autocorrelation coefficients of maximally flat filters, it
is not clear how to use this description to obtain other designs. For
example, when obtaining coiflets, neither her method
\cite{DAUBEC:1993} nor the simpler description in \cite{MO-BE-HE:1999}
can prevent bad behavior of the frequency response or guarantee zeroes
in locations different than $\pi.$ Theoretical results are also
hard to obtain, and we cannot even establish the existence of coiflets
for an arbitrary number of vanishing moments.

Here we propose a different approach for the design and implementation
of \qs. We first isolate a specific property of interest, like
interpolation
or linear
phase, and then parameterize all possible rational \qs\ with that
property. The parameters can be used to impose additional
properties or to optimize a particular design
\cite{MONZON:1994}.

In this paper we  choose linear phase
as the initial property in our design.  It implies two
possible symmetries for the wavelet and scaling functions. Either the
scaling function is even (and the wavelet is even about $1/2$) or it
is even about $1/2$ (and the wavelet is odd about $1/2$). These two
situations give rise to two different classes of linear phase filters
that we denote 0-SYM and 1/2-SYM.  In this paper
we discuss the 0-SYM case. For the 1/2-SYM case
see \cite{HER-VET:1993, SELESN:1999, MONZON:1994}.
Even though there are no linear
phase interpolating filters, it is possible to obtain linear phase and
almost interpolating filters. It is enough to ask for a symmetric
wavelet with enough vanishing moments. Both properties imply
that the scaling function has vanishing moments (twice the number for
the wavelet) and that assures the almost interpolating property
\cite{MO-BE-HE:1999}.

Since linear phase filters have poles inside and outside the unit
circle, we indicate how to implement these IIR by approximating them
by cascades of FIR filters as in~\cite{BEYLKI:1995}. The FIR
approximation does satisfy the quadrature mirror condition with
any desired accuracy. The factors in the approximation directly depend
on the location of the poles of the IIR filter, what allow one to
introduce the accuracy as part of the filter design process.

The main result of this paper is the complete description, with closed
form expressions, of all possible 0-SYM filters generating wavelets
and scaling functions with symmetry properties and any prescribed
number of vanishing moments. Each 0-SYM filter is described using a
minimal set of parameters which are directly related with the zeroes
of the filter.  This property allow us to construct filters with any
preassigned set of zeroes in the stopband. In particular, we
construct, for each fixed order, the filters which generate wavelets
with maximal number of vanishing moments. We also give sufficient
conditions on the parameters to assure that the frequency responses of
the filters are always positive or that the filters have all their
poles on the imaginary axis.

The paper is organized as follows.  In the next section we 
state some basic results and definitions and
then identify the two kinds of symmetry.
In Section~\ref{o-sym} we describe, in closed form, all solutions of the
0-SYM case. We then identify a minimal set of parameters
defining each solution, use them to obtain a variety
of properties, and show explicit examples.  
The properties of the maximally flat 0-SYM filters and  an
example of its  FIR implementation are discussed in Section \ref{MFF}.

\Section{Preliminaries}
\lb{prelim}

\begin{itemize}

\item Unless otherwise indicated, $x$ and $\xi$  are real variables while 
$z$ is a complex variable.

\item A \q\ is a $2\pi$-periodic 
function $m_{0},$
\beq{qmf0}
m_{0}(\xi) =  \sum_{k\in \mathbf{Z}} h_{k} {\rm e}^{-i k \xi}, 
\eeq
such that
\beq{qmf1}
| m_0(\xi)|^2 + | m_0(\xi+ \pi)|^2 = 1. 
\eeq
The numbers $\{ h_{k} \}$ are the \emph{coefficients} of the filter
$m_{0}.$ We assume $m_0$ to be a rational real trigonometric function.
That is, $\{ h_{k} \}$ are real and have exponential decay, but maybe
an infinite number of them is non-zero.  Condition \rf{qmf1} is then
equivalent to
\beq{qmf2}
2 \sum_{k} h_k h_{k+2n} = \delta_{no}  \bt{for} n \in \mathbf{Z},
\eeq
where $\delta_{nm}$ is defined as $\delta_{nm} = 1$
if $m = n,$ and $\delta_{nm} = 0 $ otherwise.

We denote by $H$ the z-transform of $\{ h_{k} \},$ 
$H(z) = \sum_k h_k z^k.$ We also refer to such $H$ as 
a \q\@.  Its \emph{frequency
response} is $H(e^{-i \xi}) = m_0(\xi).$

From~\rf{qmf2}, $H$ satisfies the functional equation (the \q\ equation),
\beq{QMF}
H(z)H(z^{-1})+H(-z)H(-z^{-1})=1 .
\eeq

In order to  generate a regular Multiresolution Analysis (see
\cite{COHEN:1990}), we need two additional
properties for $H.$

The first one is the \emph{normalization} or \emph{low-pass} condition. 
It forces 
\beq{qmf3}
m_{0}(0) = 1  \bt{for} H(1) = 1.
\eeq
The second one assures
that $H$ is non-zero in certain locations of the unit circle~\cite{COH-SUN:1993}.
\vs
\noindent
{\bf Cohen's condition}
There exist no non-trivial invariant cycles for $\xi \rightarrow 2\xi$
such that $ |m_0| = 1 $ on the whole cycle, i.~e., there is no $\w\neq 0$
with $2^{n}\w=\w\pmod{2\pi}$ such that for all $0\leq j < n$
\beq{qmf4}
m_{0}(2^{j}\w + \pi )  = 0.
\eeq
In practice, we first find a normalized $H$ satisfying the \q\ equation,
and then verify Cohen's condition.

\item
A solution $\varphi$ of
\beq{s1}
\varphi (\frac{x}{2}) = 2\sum_{k} h_{k} \varphi (x-k)
\eeq
is called a ``scaling function''.  Equivalently, 
\beq{s2}
{\hat \varphi} (2\xi)= m_{0}(\xi){\hat \varphi} (\xi),
\eeq
where ${\hat \varphi}(\xi) = \int_{-\infty}^{+\infty} 
\varphi(x){\rm e}^{-i \xi x}dx$ and $m_{0}$ is the \emph{associated} filter.

From~\rf{qmf3} and~\rf{s2}, 
\beq{s3}
{\hat \varphi} (\xi) = \prod_{k=1}^{\infty}
\ m_0(2^{-k}\xi)
\bt{and} 
{\hat \varphi} (0) = 1 .
\eeq

The {\em mother}
wavelet $\psi$ is defined on the Fourier side as 
\beq{uni4}
{\hat \psi} (2\xi)= m_{1}(\xi){\hat \vp} (\xi)
\eeq
where $m_1$ is a $2\pi$ periodic function, 
\beq{uni5}
m_1(\xi) = {\rm e}^{-{\rm i}  \xi} \lambda(\xi) \overline{m_0 (\xi
+ \pi)}.
\eeq
Here $\lambda$ is $\pi$-periodic and $|\lambda(\xi)|=1$ a.e. and
$\overline{m_0}$ is the complex conjugate of $m_0$.

\end{itemize}

\subsection{The moment condition}
One fundamental property for wavelet bases ~\cite{BCR, MEYER:1992,
MO-BE-HE:1999} is the property of vanishing moments of the wavelet
or the scaling function:
\beqa
&& \int_{\bf R} x^{k} \psi(x) \, dx = 0 \bt{for} 0\leq k <
M, \lb{mc1} \\
&& \int_{\bf R} x^{k}\vp(x)\, dx
=\delta_{k0} \bt{for} 0\leq k < N. \lb{mc4}
\eeqa
In terms of the filter they become \cite[Page 258]{DAUBEC:1992}
\beqa
&& D^{k}H(-1) = 0  \bt{for} 0\leq k < M, \lb{mc3} \\
&& D^{k}H(1) = \delta_{k0} \bt{for} 0\leq k < N. \lb{mc6}
\eeqa
\subsection{The symmetry condition}

We now find which filters  generate symmetry
conditions on the scaling and wavelet functions. 
For the scaling function to be symmetric, 
we need real constants $\lambda$ and
$c$ such that for all $x$ in {\bf R} $\vp (\lambda - x) =  c\,\vp
(\lambda + x),$ 
that is
\beq{sp.0}
{\hat \vp} (-\xi) = c\,{\rm e}^{2{\rm i} \lambda \xi}{\hat \vp} (\xi).
\eeq
By~\rf{s3} $c$ should be $1$ 
while using~\rf{s2} and~\rf{sp.0}
with $\xi$ and $2\xi,$
it follows  that symmetry for $\vp$ is equivalent to
\beq{sp.2}
m_{0}(-\xi) = {\rm e}^{2{\rm i} \lambda \xi} m_{0}(\xi)  .
\eeq

Since $m_{0}$ has real coefficients, 
$ m_{0}(-\xi) = \overline{m_{0}(\xi)}$ and we can write
$m_{0}(\xi ) =  a(\xi) {\rm e}^{-{\rm i} \lambda \xi},$
where $a(\xi) = {\rm e}^{{\rm i} \lambda \xi} m_{0}(\xi)$ takes real
values. In conclusion, symmetry for $\vp$ is equivalent with
linear phase for $m_{0}.$

In terms of $H$,  \rf{sp.2} becomes
\beq{sp.3}
H(z)=z^{2\lambda}H(z^{-1}).
\eeq
Since replacing  $\vp$ by an integer translate still
generates the same
Multiresolution Analysis, we can consider $\lambda$ in the interval
$[0,1),$  
and then the only possible values for $\lambda$ are $0$ or $1/2$.

Thus, either the
scaling function is even (0-symmetry) or it is even about$\;1/2$
(1/2-symmetry).  By~\rf{QMF} and~\rf{sp.3} the corresponding
conditions on the \q\ are

\begin{description} 
\item[0-SYM]
\beqa & &  H^{2}(z) + H^{2}(-z) = 1 \lb{sp.5} \\
& & H(z)=H(z^{-1}). \lb{sp.6} 
\eeqa

\item[1/2-SYM] 
\beqa
& & H^{2}(z) - H^{2}(-z) = z \lb{sp.7} \\
& & H(z)=zH(z^{-1}).\lb{sp.8}
\eeqa
\end{description}

We now show which kind of symmetry we
impose  on the
wavelet if we start with a scaling function with symmetry 
$\lambda = s_{\vp}.$ 

As in~\rf{sp.0}, we require
\beq{sp.9} 
{\hat \psi} (-\xi) = c\,{\rm e}^{2{\rm i} s_{\psi} \xi}{\hat \psi} (\xi),
\eeq
and then $c$ is $\pm 1$  because $\psi$ has
${\bf L}^2$ norm equal to $1$.

Using~\rf{uni4} and~\rf{uni5},
\[{\hat \psi} (2\xi)=   \overline{m_{0}(\xi + \pi )}
{\rm e}^{-{\rm i} \xi} d(2\xi) {\hat \vp} (\xi) \]
where $d(\xi)$ is a rational
trigonometric function and $d(\xi)d(-\xi)=1$.

Because of~\rf{sp.0} and~\rf{sp.2}
\[
{\hat \psi} (-2\xi) =  \frac{{\rm e}^{{\rm i} \xi}
\overline{m_{0}(-\xi -\pi )}  {\hat \vp} (-\xi)}{d(2\xi)} 
 = \frac{{\rm e}^{2{\rm i} \xi} {\rm e}^{-2 \pi {\rm i} s_{\vp}} 
{\hat \psi} (2\xi)}{d^{2}(2\xi)}.
\]

Hence, to obtain~\rf{sp.9} we need $d(\xi)= \pm {\rm e}^{{\rm i} m \xi} $
with $m \in {\bf Z}$ and
$ 2s_{\psi} = 1-2m $ .
Choosing $d(\xi) = -1$, that is 
\beq{m1}
m_{1}(\xi) = -{\rm e}^{-{\rm i} \xi}\overline{m_0 (\xi + \pi)},
\eeq
we obtain $s_{\psi} = 1/2$ , unique  up to integer translations. Consequently,
$\psi$ is even about $1/2$ or odd about $1/2$ depending on whether 
$s_{\vp}$ is $0$ or $1/2$. An example of the latter case is the Haar basis,
where the  wavelet and scaling
function are $\psi = 1_{[0, \half{2})} - 1_{[\half{2}, 1)}$ and 
$\vp = 1_{[0, 1)}.$
Here $1_{\rm I}$ is the indicator function of the interval ${\rm I}.$

The choice of $m_{1}$ in~\rf{m1} can be written in terms of the $z$-transforms
of the filters as
\beq{GH}
G(z)= -z H(-z^{-1}).
\eeq

\Section{The characterization of the 0-SYM case}
\lb{o-sym}

First, we introduce some additional notations and definitions.

Sometimes we  omit the variables of a function
as in the following transformations
\beqa 
\tilde{f} & = & \tilde{f}(z) = f(-z) \lb{tilde} \\
\is{f} & = & \is{f}(z) = f(z^{-1}). \lb{breve}
\eeqa

We use the notation $\sv{f}{t}$ for the sign of the function $f$ at
the value $t,$
\beq{c.1}
\sv{f}{t} = \rm{sign}(f(t))
\eeq

If $E$ is a polynomial of degree $d$
its {\em reciprocal} is $z^{d} E(z^{-1}) $.

We call $a$ an {\em all-pass} function if 
\beq{c.2}
a(z)a(z^{-1})=1.
\eeq
Rational all-pass functions can be described as: 
\beq{c.3}
a(z)= \pm z^{n}\frac{E(z)}{E(z^{-1})}
\eeq
where $n$ is an integer, $E$ is a polynomial coprime with its
reciprocal and $E$ does not vanish at zero. 

The {\em order} of a rational function $\frac{P}{Q}$, where
$P$ and $Q$ are 
coprime polynomials, is the
maximum of their degrees.

We use $\circ$ for composition of functions. If $n$ is a
positive integer, $f^{[n]}$ denotes $f\circ\ldots\circ f$
(n times), $f^{[0]}$ is the identity and $f^{[-1]}$ is the inverse of $f$. 

The bilinear transformation 
\beq{be1}
\beta(z)= \frac{1-z}{1+z}
\eeq
 maps the unit disk onto the right half
plane $\{ z: \Re(z) > 0\}$ and the unit circle onto the imaginary axis:
\beq{be2}
\beta({\rm e}^{{\rm i} \xi}) = -{\rm i} \tan (\frac{\xi}{2}).
\eeq

The transformation $\beta$ gives a natural way
to relate analog with digital filter design; see \cite{STEIGL:1965} and
\cite[Page 219]{RAB-GOL:1975}, for more details. 
The conjugation by $\beta$ of $f$ is
$ f_{\beta} = \beta\circ f \circ\beta^{[-1]} $ and similarly for
$\tilde{\beta}(z) = \beta (-z)$.

The next theorem provides an explicit description of all 0-SYM
filters. Conjugation by $\beta$ is the key element for the simplicity
of this parameterization.  In Corollary \ref{coro1} below we explain
how to obtain vanishing moments for the associated wavelet and scaling
functions.

\begin{thm}
\lb{thm1}
There is a bijection between the set of 0-SYM filters $H$ and the set of
even all-pass functions $a$ given by:
\beqa
a & \longmapsto & H = 
\beta\circ\frac{1}{(1+\sqrt{2}\, a)^{2}}\circ\beta^{[-1]}, \lb{Ha} \\
H & \longmapsto & a = 
\frac{1+H(z)-H(-z)}{\sqrt{2}H(-z)}\circ\beta . \lb{aH}
\eeqa
If $H(1)=1$, $H$
can be described as
\beq{HA}
H(z) =  \frac{A(z)(A(z)\pm \sqrt{2}A(-z))}
{A(z)^{2}+A(-z)^{2}\pm \sqrt{2}A(z)A(-z)},
\eeq
where $A(z)$ is a polynomial of even degree, coprime
with $A(-z),$ and equal to its reciprocal. Also,
 $A(0) A(1) \neq 0 $ and $A(-1) = 0.$
\end{thm}

\noindent {\bf Proof}
The transformation $\beta$ satisfies:
\beqa
&& \beta^{[-1]} = \beta, \lb{be7} \\
&& \is{\beta} = -\beta, \lb{c.13} \\
&& \tilde{\beta} = \inv{\beta}.\lb{c.14}
\eeqa
Conjugation by $\beta$ yields,
\beqa
d \;\; \mbox{is all-pass} & \Longleftrightarrow & d_{\beta} \;\; \mbox{is
an odd function},  \lb{c.10} \\
d = \is{d} & \Longleftrightarrow & d_{\beta} \;\; \mbox{is an
even function.} \lb{c.11}
\eeqa 
Define the function $\nu$ as
\beq{nu}
\nu(z) =  \beta\circ\frac{1}{(1+\sqrt{2}\, z^{-1})^{2}} =
 \frac{1 + \sqrt{2}z}{1 + \sqrt{2}z + z^{2}}.
\eeq
First we  show that the maps in \rf{Ha} and \rf{aH} are well defined.
We need to check
that for $a$ even all-pass 
\beq{Hnu}
H = \nu\circ a(z^{-1})\circ\beta = \nu\circ\beta\circ (\frac{1}{a})_{\beta}
\eeq
belongs to 0-SYM.
Because $a$ is even, $a_{\beta}$ equals
$(a_\beta)\,\is{}$, 
which gives \rf{sp.6}. On the other hand, ~\rf{sp.5} is invariant
under 
odd functions and
$(\frac{1}{a})_{\beta}$  is odd, hence 
we only need to verify that $\nu\circ\beta$ satisfies \rf{sp.5}.
This follows  replacing $z$ by $\beta(z)$ in the identity
\beq{nu2}
\nu^{2}(z) + \nu^{2}(z^{-1}) = 1.
\eeq

Next we start with  $H$ with 0-SYM and check
that $a$ in \rf{aH} is even and all-pass.
Using~\rf{c.10} and~\rf{c.11} it suffices to show that
\[ \beta\circ\frac{1+H-\tilde{H}}{\sqrt{2}\tilde{H}} \]
is odd and invariant under $z^{-1}$.
The latter follows from \rf{sp.6}, while using \rf{c.13} the former
is equivalent to the identity
\[ \beta\circ\frac{1+\tilde{H}-H}{\sqrt{2}H}=\beta\circ
\frac{\sqrt{2}\tilde{H}}{1+H-\tilde{H}} \]
or
\[ 1 - (\tilde{H}-H)^{2} = 2H\tilde{H}  \]
which follows from \rf{sp.5}. 

To establish that the maps in \rf{Ha} and \rf{aH} are each others
inverses, let us denote them by ${\cal S}$ and ${\cal T}$
respectively.

Let $a$ be an even all-pass function. Since
\[ {\cal S}(a) = \nu\circ \is{a}\circ\beta \]
we have
\[ \frac{1+{\cal S}(a)-{\cal S}(a)(-z)}{\sqrt{2}{\cal S}(a)(-z)} =
(\frac{1+\nu (\is{a}) - \nu (a)}{\sqrt{2}\nu (a)})\circ\beta . \]
Consequently, using that $a$ is all-pass and the
identity
\[ \frac{1+\is{\nu}  - \nu }{\sqrt{2}\nu } = z , \] 
we have
\[ {\cal T}\circ {\cal S}(a) =
(\frac{1+\is{\nu} - \nu }{\sqrt{2}\nu })(a) = a. \]

Finally, if $H$ is a  0-SYM filter,
\[
1+\sqrt{2}{\cal T}(H) =  (\frac{\tilde{H}}{1-H})\circ\beta, \]
\beqan
{\cal S}\circ {\cal T}(H) 
& = & \beta\circ (\frac{(1-H)^{2}}{\tilde{H^{2}}}) =
\frac{\tilde{H^{2}}-(1-H)^{2}}{\tilde{H^{2}}+(1-H)^{2}} \\
& = & \frac{1-H^{2} - (1+H^{2}-2H)}{1-2H +(H^{2}+\tilde{H^{2}})} = H.
\eeqan

For the second part of the theorem, 
we write the even all-pass function $a$ as

\beq{aB}
a(z)= \sv{a}{i} (-1)^{m + r} \frac{B(z^2)}{z^{2m+2r}B(z^{-2})}
\eeq
where $r \geq 0$, $m \geq 1$ (because $H(1)=1$),
$B(z)=\sum_{j=0}^{r} b_{j}z^{j}$ with $b_{0}b_{r}\neq 0,$
and $B$ coprime with its reciprocal. In particular 
$B(1)\neq 0$.

Thus,
\beq{dt.1}
\is{a}(\beta(z)) = \sv{a}{i} (-1)^{m + r} \frac{A(-z)}{A(z)} 
\eeq
where
\beqa
 A(z) & = & (1+z)^{2m+2r}B(\beta^{2}(z)) \lb{dt.3} \\
& = & (1+z)^{2m} \sum_{j=0}^{r} b_{j}(1-z)^{2j}(1+z)^{2r-2j} .
\lb{dt.2}
\eeqa
The properties for $A$ follow from the conditions on $B$ while~\rf{HA}
is the result of applying $\nu$ to~\rf{dt.1}. \ep

As stated in the introduction, we can obtain an almost interpolating
0-SYM filter by 
simply choosing an associated all-pass function with a high
number of poles at $0.$

\begin{corollary}
\lb{coro1}
Consider $a$ and $H$ as in the previous
theorem and let $m>0.$  The following
conditions are equivalent:

\begin{description}
\item[1]  The multiplicity of $0$
as a pole of $a$ is $2m.$
\item[2] The multiplicity of $-1$
as  a zero of $H$ is $2m.$
\item[3] The multiplicity of $1$
as  a zero of $1-H$ is $4m.$
\end{description}

If any of these conditions holds then by \rfE{mc3} and~\rfE{mc6}
\beqan
&& \int_{\bf R} x^{k} \psi(x) \, dx = 0 \bt{for} 0\leq k < 2m,
\bt{and} \\
&&
\int_{\bf R} x^{k}\vp(x)\, dx
=\delta_{k0} \bt{for} 0\leq k < 4m,
\eeqan
where $\psi$ and $\vp$ are the wavelet and scaling functions associated
to $H$.
\end{corollary}
\noindent {\bf Proof }
We simply show that the three conditions are equivalent
with the polynomial $A$ of~\rf{dt.2} having $2m$
zeroes at $-1.$ Part \textbf{1} follows from \rf{aB} and \rf{dt.2}
while~\rf{HA} gives Part \textbf{2} (note that $A(1) \neq 0 $ but
$A(1) = 0 $ because $m > 0$).
The equivalence with Part \textbf{3} holds because~\rf{HA}
yields
\beq{1-H}
1-H(z)= \frac{A(-z)^{2}}{A(z)^{2}+A(-z)^{2}\pm \sqrt{2}A(z)A(-z)}
\eeq
where the polynomials in the numerator and denominator are coprime. \ep


\begin{rems}
\lb{rem1}
{\em

\begin{enumerate}
\item
We can use the  theorem as a dictionary to reformulate
properties of $H$ in
terms of $a$ and vice versa.  Let us  list some of them:
\beqa 
z_{0} \in H^{-1}\{1\} & \Longleftrightarrow & 
a(\beta (z_{0})) = \infty \lb{H+} \\
z_{0} \in H^{-1}\{-1\} & \Longleftrightarrow & a(\beta (z_{0})) =
 -\frac{1}{\sqrt{2}} \lb{H-}   \\
z_{0} \in H^{-1} \{ \infty \} & \Longleftrightarrow & a(\beta (z_{0})) =
{\rm e}^{\pm {\rm i} \pi 5/4 } \lb{H++} 
\eeqa
Because of symmetry and the \q\ condition, the frequency response of
$H$ only takes values in the real interval $[-1, 1].$ In
Table~\ref{tab1} we listed the correspondence between values of $H$ on
the unit circle and values of $a$ on the imaginary axis. This
correspondence can be easily obtained from~\rf{Hnu}.

\item Because of \rf{sp.5}, the zeroes of $H$ are
the negatives of the preimages by $H$ of $1$ and $-1.$
These two sets of preimages are nicely represented in the numerator
of~\rf{HA}. The zeroes of $A(-z)$ (counted twice) are exactly the set
$H^{-1}\{1\}$ (see \rf{1-H}) 
while the set $H^{-1}\{-1\}$ is given
by the zeroes of $A(-z) \pm \sqrt{2} A(z)$ (counted twice) because
\beq{1+H}
1 + H(z) = 
\frac{(A(-z) \pm \sqrt{2} A(z))^{2}}{A(z)^{2}+A(-z)^{2}\pm \sqrt{2}A(z)A(-z)}.
\eeq
The even multiplicity of the preimages of $1$ and $-1$ agrees with the
fact that the frequency response of $H$ cannot take values outside
$[-1, 1].$

\item The filter $H$ cannot have poles either
on the unit circle or on the real line, because $\beta$ maps those values
to the imaginary and real axes, and $a$, being even and real, maps them 
to the real line.  

Also,  if $H$ were a polynomial,
$H^{-1}\{\infty\} = \{\infty\}$, and we would have
\[ a(-1) = {\rm e}^{\pm {\rm i} \pi 5/4  }, \]
leading to a contradiction because $a$ has real
coefficients. Therefore, we recover the result on the absence of real
polynomial \qs\ generating wavelets with even
symmetry~\cite{SMI-BAR:1986, DAUBEC:1988}.
\item Using the notation in \rf{c.1}, it follows from \rf{aB}
that $\sv{a}{1} = (-1)^{m + r} \sv{a}{i}.$

Since $H(0)$ and $H(i)$ are real numbers, \rf{sp.5} implies
\beq{r.2}
H(0) = \frac{\sv{H}{0}}{\sqrt{2}} \bt{and} 
H(i) = m_{0}(\pi/2) = \frac{\sv{H}{i}}{\sqrt{2}}. 
\eeq

Finally,~\rf{aH} gives
\beq{r.4}
\sv{a}{1} = \sv{H}{0} \bt{and} \sv{a}{i} = \sv{H}{i}. 
\eeq

\end{enumerate}
} 
\end{rems}

\renewcommand{\arraystretch}{1.5}
\begin{table}[htb] \caption{\lb{tab1} Correspondence between  values
of $H$ 
and $a$}
\begin{center}
\begin{math}
\begin{array}{ r  l}\hline
H(e^{\rm{i} \xi}) & a(-{\rm i} \tan (\half{\xi}) ) 
\\
\hline
\{ -1 \} & \{ -\inv{\sqrt{2}} \} \\
(-1, 0) & (-\sqrt{2}, -\inv{\sqrt{2}})  \cup (-\inv{\sqrt{2}}, 0) \\
\{ 0 \} & \{ -\sqrt{2}, 0 \} \\
(0, 1) & (-\infty, -\sqrt{2}) \cup (0, +\infty) \\
\{ 1 \} & \{ \infty \} \\
\hline
\end{array}
\end{math}
\end{center}
\end{table}
\renewcommand{\arraystretch}{1}

We now identify a minimal set of parameters defining each 0-SYM
filter. 
\vs
\noindent {\bf Parameterization in terms of preimages of one}

Assume $H(1)=1$. From~\rf{H+}, the set of preimages 
of $1$ counted with their 
multiplicity (the set $H^{-1}\{1\}$) completely
determines the poles of the all-pass function $a$. 
Therefore, up to a sign, $a$, and thus $H$, are
completely determined by $H^{-1}\{1\}$. We fix the sign
by choosing $\sv{a}{\rm{i}},$ that is $\sv{H}{\rm{i}}.$

The conclusion is that two solutions of 0-SYM are the same if their sets 
of preimages of $1$ and their values at $\rm{i}$ are identical.

We now study the degrees of freedom in the set $H^{-1}\{1\}.$

First, because of \rf{1-H}, $H^{-1}\{1\}$ is the set of zeroes of
$A(-z)^2.$ Using the properties of $A$ in  \rf{HA} and~\rf{dt.2}, 
$A(-z)$ has $2 (m + r)$ zeroes,  $2 m$ of them
being
$1.$ Also, if $w$ belongs to $H^{-1}\{1\}$, $\overline{w}$ and $w^{-1}$
belong
to $H^{-1}\{1\}$ but $-w$ does not. 

Therefore, once $\sv{H}{i}$ and $m$ are fixed, there are only $r$
degrees of freedom for a 0-SYM filter of order $4 (m + r).$
We isolate these $r$ parameters in  a
subset $\Lambda$ of $H^{-1}\{1\}$ by
discarding reciprocals and the value $1,$
\beq{G}
\Lambda = \{ \lm_1, \cdots, \lm_r  \}.
\eeq

In the next theorem, we rewrite our previous parameterization
of $H$ to show the dependence on $\Lambda.$

\begin{thm}
\lb{thm2}
\vs
Let $H$ be any 0-SYM filter with $2 m$ zeroes at $-1$ and
$\Lambda$ the subset of $H^{-1}\{1\}$ defined in~{\rm \rf{G}}.
If $\nu$ is the function defined in~{\rm \rf{nu} } and $\eta$ is
the function
\beq{eta}
\eta(z) =  \frac{z + z^{-1}}{2},
\eeq
then
\beq{H1}
H(z) = \nu\left( \sv{H}{i} \beta(\eta(z))^m \prod_{j = 1}^r
\beta(\frac{\eta(z)}{\eta(\lm_j)}) \right).
\eeq
\end{thm}
\noindent {\bf Proof }
Since \rf{Hnu} and \rf{dt.1} imply
\[
H(z) = \nu(a(\frac{1}{\beta})) 
= \nu(\sv{a}{i} (-1)^{m + r} \frac{A(-z)}{A(z)}) , \]
it suffices to show that
\beq{ap}
\sv{a}{i} (-1)^{m + r} \frac{A(-z)}{A(z)} = u(\eta(z))
\eeq
where
\beq{u}
u(z) = \sv{a}{i} \beta(z)^m \prod_{j = 1}^r
\beta(\frac{z}{\eta(\lm_j)}).
\eeq	

We know that $A$ has degree
$2 (m + r),$ it is equal to its reciprocal, and (from \rf{1-H}) its
zeroes 
are the negatives
of the elements of the set $\Lambda$ considered together with their
inverses.
Therefore,  there exist a constant $c$ and a polynomial $E,$
\[ E(z) = \prod_{j = 1}^r (z + \eta(\lm_j)), \]
such that
\[
A(z) = c z^{m+r} (1 + \eta(z) )^m E(\eta(z))
\]
and the theorem follows. \ep


We have now a straightforward way of designing a 0-SYM filter
$H$ by imposing the locations
 where the filter takes the value one. 

Since the value $m_{0}'(\half{\pi})$ gives an idea of how narrow the
transition band is,  we will express it in terms of $H^{-1}\{1\}.$
The next corollary implies that if we trade vanishing moments for any
other value on the arc $(\half{\pi}, \pi)$ of the
unit circle (the stopband), then the slope of $m_{0}$ is steeper at $\half{\pi}.$

\begin{corollary}
\lb{coro2}
Let $H$ as in Theorem {\rm \ref{thm2}} and $m_0(\xi) = H(e^{-i \xi}).$
Then
\beq{sl}
m_0'(\half{\pi}) = -(2 - \sqrt{2} \sv{H}{i}) (m + \sum_{j=1}^r
\inv{\eta(\lm_j)}).
\eeq

\end{corollary}
\noindent {\bf Proof }
From the proof of Theorem \ref{thm2}, we have
\beq{mu}
m_{0}(\xi) = \nu(u(\cos(\xi))).
\eeq
Taking derivatives at $\pi$ and using $u(0) = \sv{H}{I}$
we obtain
\[
m_{0}'(\half{\pi}) = -\nu'(\sv{H}{I}) u'(0) = -\half{\sqrt{2} - 2 \sv{H}{I}}
u'(0).
\]
Since
\[ \frac{\beta'(0)}{\beta(0)} = -2 \]
the logarithmic derivative of \rf{u} at $0$ is,
\beqan
\frac{u'(0)}{u(0)} = -2 (m + \sum_{j=1}^{r}
\inv{\eta(\lambda_j)}).  \hs \hs \hs \ep
\eeqan
\vs

The next two results are obtained imposing
conditions on the set $H^{-1}\{1\}.$

\begin{prop}
\lb{prop1}
Let $H$ be a 0-SYM filter with $\sv{H}{i} = 1$ and\\
$H^{-1}\{1\} \cap \{ z: |z| = 1 \} = \{1\}.$
Then
\[ H(e^{i\xi}) > 0 \bt{for} \xi \in (-\pi, \pi).
\]
\end{prop}

\noindent \textbf{Proof} $\;$
Since $H(1) = 1$ and $H(z^{-1}) = H(z)$ it suffices to  consider
$\xi$ in $(0, \pi).$ From~\rf{aB},
\beq{aE}
a(\beta(\e{\xi})) = a(-i \tan(\half{\xi})) = \frac{E(x)}{x^{m + r}
E(x^{-1})},
\eeq
where $x = \tan^2(\half{\xi}),$ $m \geq 1,$   and
$E(x) = B(-x).$ 

The zeroes of $E$ are $-\beta(-\lambda)^2$ where $\lambda \neq 1$
belongs to  $H^{-1}\{1\}.$ (See \rf{dt.3}  and the second remark in
Remarks~\ref{rem1}.)
These zeroes cannot be positive numbers because
\[
-\beta(-\lambda)^2 \in (0, +\infty) \iff \beta(-\lambda) \in
\mathbf{R} i \iff |\lambda| = 1,
\]
which contradicts our assumption on $H^{-1}\{1\}.$

Thus, $E$ has constant sign in $(0, +\infty)$ and \rf{aE} implies
that $a(\beta(\e{\xi}))$ is positive for $\xi \in (0, \pi).$ 
The result follows using Table \ref{tab1}. \ep
\vs
We now discuss a typical filter design problem.
Suppose the desired frequency response vanishes
at certain specific locations of the stopband. 
By choosing  appropriate preimages of $1$ in the
passband, we can obtain the given zeroes 
but, in principle,  we cannot avoid
further zeroes in the stopband because they are related
to preimages of $-1.$ The next theorem provides
sufficient conditions to prevent this last situation from
happening and,
incidentally, forces all poles of the filter to be purely
imaginary.
By restricting $H^{-1}\{1\}$ to the right half  plane
we can parameterize 0-SYM filters with any
preassigned set of zeroes in the stopband. We can
then obtain a filter with minimal  order or, considering
more parameters, even impose additional properties.

\begin{thm}
\lb{thm3}
\vs
Let $H$ be any 0-SYM filter with $H^{-1}\{1\} \subseteq \{ z : \Re(z) > 0 \}.$ 
We have the following properties:
\begin{description}
\item[A] $H^{-1}\{\infty\} \subseteq \mathbf{R} i,$

\item[B] $H^{-1}\{-1\} \subseteq \{ z : \Re(z) < 0 \}$.
\end{description}
\end{thm}
\noindent {\bf Proof }
The assumption on $H$ and \rf{H+} imply that all poles
of the corresponding all-pass function $a$ are in the unit
disk. Therefore, $a(z^{-1}) = \inv{a(z)}$ is a finite Blaschke product
and then $a$ maps the unit circle into itself
and $\{ z : |z| > 1 \}$ into the unit disk and vice versa.

To verify Part \textbf{A} consider $p \in H^{-1}\{\infty\}.$
By \rf{H++}, $a(\beta(p))$ and then $\beta(p)$ belong to the unit
circle. Thus, $p$ is purely imaginary by \rf{be2}.

For Part \textbf{B} start with  $p \in H^{-1}\{-1\}.$ Since \rf{H-} implies
that $a(\beta(p))$ belongs to the unit disk, we have $|\beta(p)| > 1$
and then $\Re(p) < 0.$ \ep

\begin{corollary}
\lb{coro3}
For $ \{ \theta_1 \cdots \theta_s \}$  in $(\half{\pi}, \pi)$ let 
\beq{Hz}
H(z) = \nu\left( \sv{H}{i} \beta(\eta(z))^m 
\prod_{k = 1}^t \beta(-\frac{\eta(z)}{\cos(\theta_k)})
\prod_{j = 1}^s \beta(\frac{\eta(z)}{\eta(\lm_j)}) 
 \right),
\eeq
where the set $\{ \lm_1, \cdots,
\lm_t \}$ is contained in
$\{ z: \Re(z) > 0 \; \mbox{and} \; |z| \neq 1 \} $ and satisfies
the conditions listed for \emph{\rf{G}}.

Then, $H$ is a 0-SYM filter that vanishes in the stopband 
only at the $2 (m + s) $ values 
\[
\left\{ \e{\theta_1},
e^{-i \theta_1}, \cdots, \e{\theta_s}, 
e^{-i \theta_s}, \underbrace{-1, \cdots, -1}_{2 m} \right\} .
\]
\end{corollary}

\noindent \textbf{Proof} $\;$
Using Theorem \ref{thm2}, define $H$ with $2 m$ zeroes at $-1$ and
$\Lambda = \{ -\e{\theta_1}, \cdots, -\e{\theta_s}, \lm_1, \cdots,
\lm_t \}.$
Theorem \ref{thm3} assures that there are no other zeroes
in the stopband. \ep

\vs

\begin{figure}[htb]
\epsfxsize = 500pt
\centerline{\epsffile{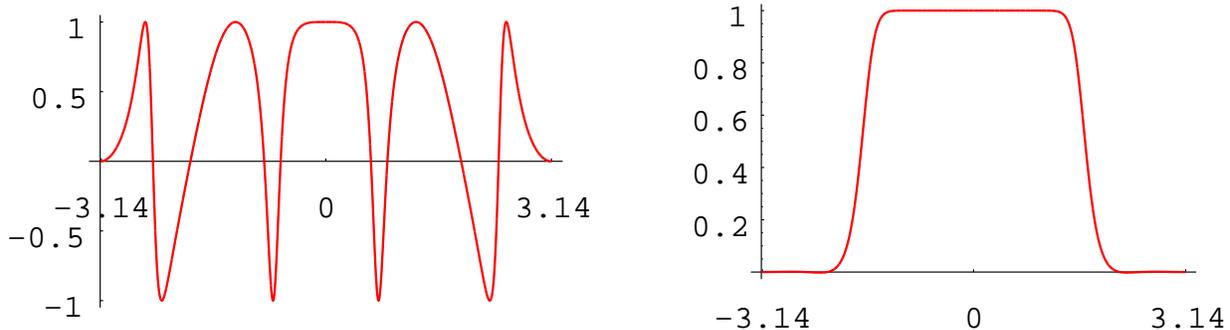}}
\caption{\label{figB} Two 0-SYM filters of order $12.$ On the left
we violated Cohen's condition. On the right we choose all preimages
of $1$ in the right half plane.}  
\end{figure}

Note that ``wrong'' choices of $H^{-1}\{ 1 \}$ can lead to very poor
frequency responses.
In Figure~\ref{figB} we plotted the frequency responses of two 
order $12$ filters with  $s^H_i =1$ and two zeroes at
$-1$ but different choices for the other zeroes.
Their set of preimages, zeroes, and poles are listed in Table
\ref{tab2}. In one  case we have chosen the zeroes in
order to violate Cohen's condition.
The poor frequency response is not a peculiarity of this choice.
Perturbing the values of the preimages we still
obtain a similar response. Note that the zeroes either 
belong to the unit circle or they are real numbers.
For the other example, all preimages belong to the right half plane 
and therefore all zeroes related to preimages of $-1$ belong
to the same half plane. In agreement with the previous
theorem, all poles are purely imaginary.

\begin{table}[htb] \caption{\lb{tab2} Zeroes and poles of filters constructed
via their preimages of one.}
\begin{tabular}{|c | l | l |} \hline
Parameters &  Zeroes  &  Poles\\ [0.5ex]
 \hline\hline
{} & -1 & -0.84955807 - 0.74802903 I \\
{} & -1 & -0.84955807 + 0.74802903 I \\
{} & -0.74212 - 0.67026705 I & -0.66304573 - 0.58380642 I \\
{} & -0.74212 + 0.67026705 I & -0.66304573 + 0.58380642 I \\
$s^H_{i} = m = 1$ & -0.30901699 - 0.95105652 I & -0.40197132 I \\
{} & -0.30901699 + 0.95105652 I & 0.40197132 I \\
{} & 0.17142917 & -2.4877396 I \\
$\Lambda = \{ e^{\frac{2}{5} \pi i }, e^{\frac{4}{5} \pi i } \}$ &
0.65396257 - 0.75652691 I &  2.4877396 I \\
{} & 0.65396257 + 0.75652691 I & 0.66304573 - 0.58380642 I \\
{} & 0.80901699 - 0.58778525 I & 0.66304573 + 0.58380642 I \\
{} & 0.80901699 + 0.58778525 I & 0.84955807 - 0.74802903 I \\
{} & 5.8333128 & 0.84955807 + 0.74802903 I \\[0.5ex]
 \hline \hline
{} & -1 & -0.083442717 I \\
{} & -1 & 0.083442717 I \\
{} & -0.79015501 - 0.61290705 I & -0.57528543 I \\
{} & -0.79015501 + 0.61290705 I & 0.57528543 I \\
$s^H_{i} = m = 1$ & -0.56208338 - 0.82708057 I & -0.73702991 I \\
{} & -0.56208338 + 0.82708057 I & 0.73702991 I \\
{} & 0.03560146 - 0.65573566 I & -1.356797 I \\
$\Lambda = \{ e^{0.21 \pi i }, e^{0.31 \pi i } \}$ & 0.03560146 + 0.65573566 I & 1.356797 I \\
{} & 0.036837087 & -1.7382676 I \\
{} & 0.082552825 - 1.5205228 I &  1.7382676 I \\
{} & 0.082552825 + 1.5205228 I & -11.984269 I \\
{} & 27.146554 & 11.984269 I \\
[0.5ex]
\hline
\end{tabular}
\end{table}
\Section{Maximally Flat Filters}
\lb{MFF}
As in Daubechies' construction for compactly
supported wavelets \cite{DAUBEC:1988}, we use our representation of 
the 0-SYM class to
construct, for each fixed order of the filter,
symmetric wavelets with maximum number of vanishing moments. For the
reasons explained before, they correspond to the choice 
$a(z) = \sv{H}{i} (-1)^n z^{-2n}$ in ~\rf{Hnu}.
For simplicity of notation, we  write $\sv{H}{\rm{i}} = (-1)^{\delta},$
where $\delta$ is a positive integer.

\subsection{Description and Properties}
\lb{MF1}

For $\delta$ and $n$ positive integers, 
we define 
\beq{E2}
E_{4n}^{\delta}(z) = \nu((-1)^{\delta + n} \beta^{2n}(z)) 
\eeq
that is

\beq{En1}
E_{4n}^{\delta}(z) = \frac{(1+z)^{2n}((1+z)^{2n}+ (-1)^{\delta + n}
\sqrt{2}(1-z)^{2n})}
{(1+z)^{4n}+(1-z)^{4n}+ (-1)^{\delta + n} \sqrt{2}(1-z^2)^{2n}}.
\eeq

The filters $E_{4n}^{\delta}$  have maximum number of zeroes at
$-1$ ($2n$ in fact) in the set of 
rational \q\ with 0-SYM and order at most $4n$.
For their frequency response, we use~\rf{be2} and~\rf{E2}:

\beqa
E_{4n}^{\delta}({\rm e}^{{\rm i} \xi}) 
& = & 
      \nu((-1)^{\delta} \tan^{2n}\frac{\xi}{2}) \lb{E2.5} \\
& = & 
      1 - \frac{1}{\half{1} + (\inv{\sqrt{2}} + (-1)^{\delta} \cot^{2n}\frac{\xi}{2})^{2}}. 
\lb{E3}
\eeqa

It readily follows that,
as $n$ goes to $\infty$,
$E_{4n}^{\delta}$  is approaching the ideal low-pass filter $\iota$

\beq{sy4}
\iota(\xi) = \left\{\ba{lll} 1 & \mbox{if} & 0\leq |\xi | < \pi/2, \\
0 & \mbox{if} & \pi/2 < |\xi | < \pi. \ea \right.
\eeq

\begin{figure}[htb]
\epsfxsize = 500pt
\centerline{\epsffile{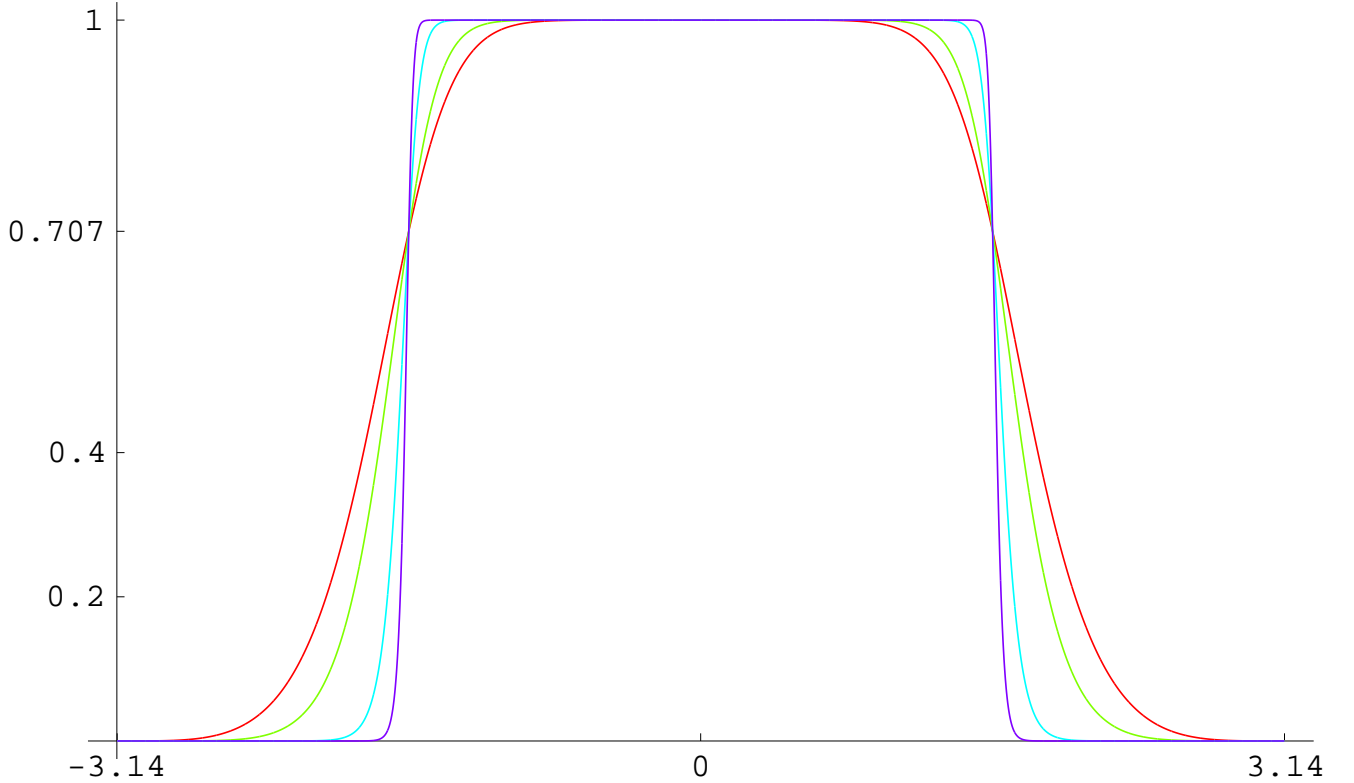}}
\caption{\label{fig1} Approximating the ideal low-pass filter by
maximally flat
filters $E^{0}_{4n}.$ The filters shown in the picture correspond to $n = 2, 3, 8,$ and $20.$}
\end{figure}

\begin{figure}[htb]
\epsfxsize = 500pt
\centerline{\epsffile{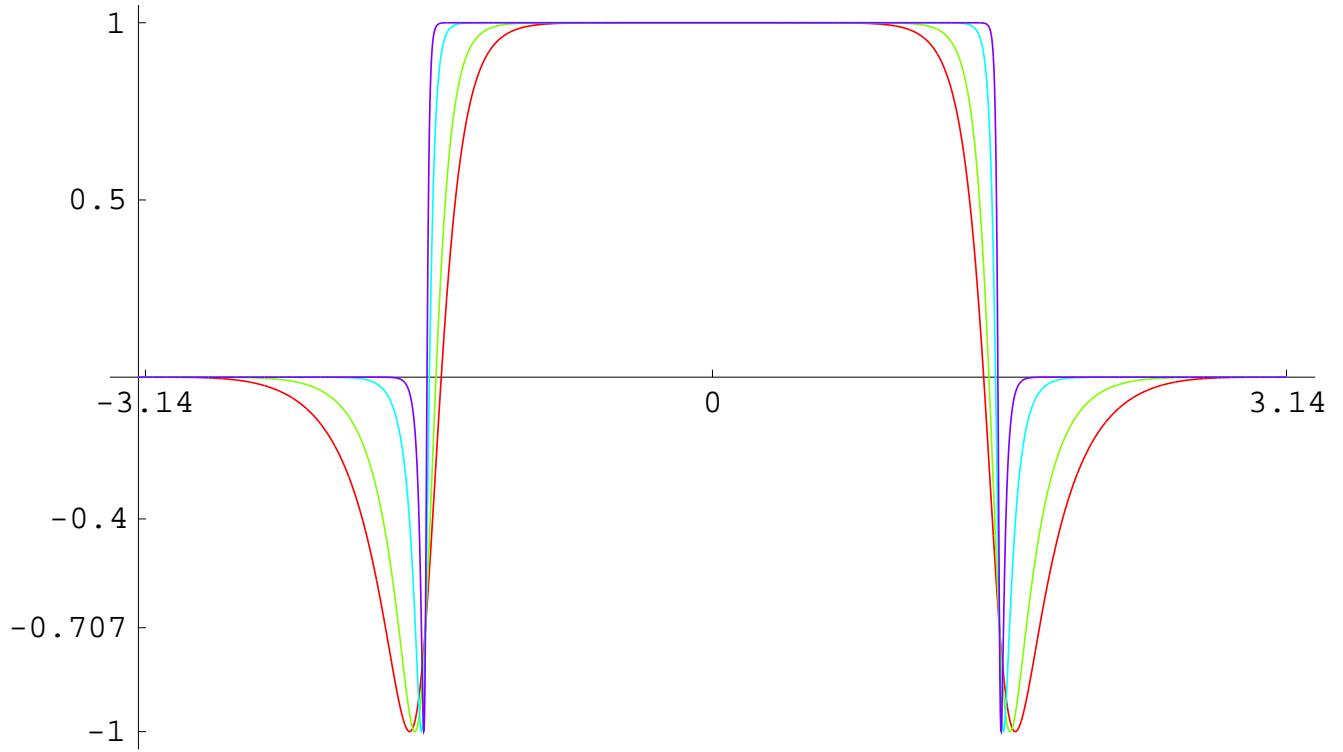}}
\caption{\label{fig2} Approximating the ideal low-pass filter by
maximally flat 
filters $E^{1}_{4n}.$ The filters shown in the picture correspond to
$n = 2, 3, 8,$ 
and $20.$}
\end{figure}

When $\delta = 0$ the
convergence is uniform because $\nu$ is a decreasing function in
the interval $(0, 1).$ To illustrate the convergence, a few filters are 
simultaneously plotted in Figure \ref{fig1} (for  $\delta=0$) and 
Figure \ref{fig2} (for  $\delta=1$).

Because of \rf{GH}, the filters $m_1$ associated with
the wavelet converge to the ideal high-pass filter. Compare this
situation with Daubechies' family of maximally flat filters that converge to the
ideal filters only in absolute value \cite{LAI:1995}.

Next, we explicitly find all zeroes and poles of $E_{4n}^{\delta}$.
In agreement with the results of the previous section, the frequency
response of $E_{4n}^{0}$ is positive in $(-\pi, \pi)$ (Proposition
\ref{prop1}), all poles of $E_{4n}^{\delta}$ are purely imaginary
(Theorem \ref{thm3}, Part \textbf{A}), and except for $-1,$ all zeroes of
$E_{4n}^{\delta}$ belong to the 
right half plane (Theorem \ref{thm3}, Part \textbf{B}). The
latter  result
also holds for Daubechies' family of maximally flat filters
\cite[Theorem 2]{SHE-STR:1996}.

\vs
\noindent {\bf Zeroes of} $E_{4n}^{\delta}$
\vs
We already know that $-1$ is a zero of multiplicity
$2 n.$ The other zeroes are
the values of $\w$ such that $H(-\w) = -1$. 
Since $a(z) = (-1)^{\delta + n} z^{-2n}$, using~\rf{H-} we need to solve
\[  \beta^{2n}(\w) = \frac{{\rm e}^{{\rm i} \pi (n + \delta - 1)}}{\sqrt{2}}. \]

It follows that  the zeroes are
\beq{E4}
\beta\left(2^{-\frac{1}{4n}}
{\rm e}^{\frac{{\rm i} \pi}{2n} (2j + n + \delta -1)}\right),
\eeq
for $-n < j \leq n.$
Note that they belong to the right half-plane.

With respect to zeroes on the unit circle, we 
now show that there is a pair of
complex conjugate zeroes if $\delta = 1$ and no zeroes if
$\delta = 0$.

We check in \rf{E4} when the argument of $\beta$ is purely imaginary.
Since $-n < j \leq n,$ it follows that $\sin(\frac{\pi}{2n}(2j+\delta
-1)) = 0$ only if $\delta = 1$ and either $j = 0$ or $j = n$.  In
either case, the real part equals
\[ \frac{\gamma}{\eta}(2^{\inv{4n}})=\beta(2^{-\frac{1}{2n}}).\]
For all $n$ the numbers above are less than $\half{1}$  and Cohen's
condition holds.  

Note that for  any real numbers $\lambda$ and $\xi$,
\beq{be5}
\beta(\frac{{\rm e}^{{\rm i} \xi}}{\lambda}) = \frac{\gamma(\lambda) - \sin\xi
 \, {\rm i}}{\eta(\lambda) + \cos\xi}
\eeq
where $\eta$ is the function in~\rf{eta} and 
$\gamma(z) =  \frac{z - z^{-1}}{2}.$

Thus, except for $-1$, the zeroes of $E_{4n}^{\delta}$ can be 
also expressed as
\beq{E5}
\left\{\frac{\gamma(2^{\inv{4n}}) -{\rm i}
\cos(\frac{\pi}{2n}(2j+\delta -1))}{\eta(2^{\inv{4n}})
- \sin(\frac{\pi}{2n}(2j+\delta -1))}\right\},
\eeq
for $-n < j \leq n.$

\noindent {\bf Poles of} $E_{4n}^{\delta}$

By~\rf{H++}, we  need the values of $z$ such that
\[ \beta^{2n}(z) = {\rm e}^{\pm {\rm i}\pi(5/4 + \delta + n)}.\]
Because of \rf{be2} and \rf{be7} the poles are
\[ 
\left\{\pm {\rm i} \tan \left (\frac{\pi}{16n}(5+8j+4\delta + 4n)\right )\right\} 
\]
for $-n <  j \leq n .$

\begin{rem}
\lb{rem2}
\em 

The recurrence relation
\[ 
E_{8s}^{\delta} = \nu((-1)^{\delta}\beta^{4s}) = E_{4s}^{\delta +
s}(\frac{1}{\eta}),
\]
follows substituting $z$ by $\inv{\eta}$ in \rf{E2}
and using the identity
\beq{dy1}
\beta\circ\frac{1}{\eta}\circ\beta^{[-1]} = z^{2}.
\eeq

Thus, for all $s$,
\[
E_{2^{n} (4 s)}^{\delta} = E_{4s}^{\delta + s}((1/\eta )^{[n]}),
\]
and the behavior of $E_{4n}^{\delta}$ for large $n$ is dictated by the
 iteration of the map $1/\eta$. This map is conjugate to $z^2$
 implying that its Julia set is the unit circle.  We have then another
 explanation why these families are approaching the ideal filter:
 When $n$ goes to $\infty$, $(1/\eta)^{[\infty]}$ restricted to the
 unit circle has only three possible values. In fact,
\beq{dy2}
(1/\eta)^{[\infty]}(e^{i \xi}) = 
\left\{\ba{lll} 1 & \mbox{if} & |\xi | < \pi/2, \\
\infty & \mbox{if} & |\xi | = \pi/2, \\
0 & \mbox{if} & \pi/2 < |\xi | \leq \pi. \ea \right.
\eeq
Theorem~\ref{thm1} implies that 
\[
H_{\beta} =  \beta\circ H \circ\beta^{[-1]} =
\frac{1}{(1+\sqrt{2}a)^{2}}
\]
for any 0-SYM filter $H.$
When $H(1) = 1$, this conjugation changes
the location of the fixed point from $1$ to
$0.$
We can use $\tilde{\beta} = \beta(-z)$ to change it from $1$ to $\infty$
\[
H_{\tilde{\beta}}(z)  = (1+\sqrt{2}a(z^{-1}))^{2} .
\]
The choice $a(z) = (-1)^{\delta + n} z^{-2n}$ yields the maximal family
$E_{4n}^{\delta}$
whose set of
preimages of $1$ equals $\{ \underbrace{1, \cdots, 1}_{4 n} \}.$
Thus,  the conjugation by
$\tilde{\beta}$ of each $E_{4n}^{\delta}$ gives a polynomial! It is
\beq{dy5}
(1+(-1)^{\delta + n}\sqrt{2}z^{2n})^{2} .
\eeq
\em  
\end{rem}
\begin{figure}[htb]
\epsfxsize=500pt
\centerline{\epsffile{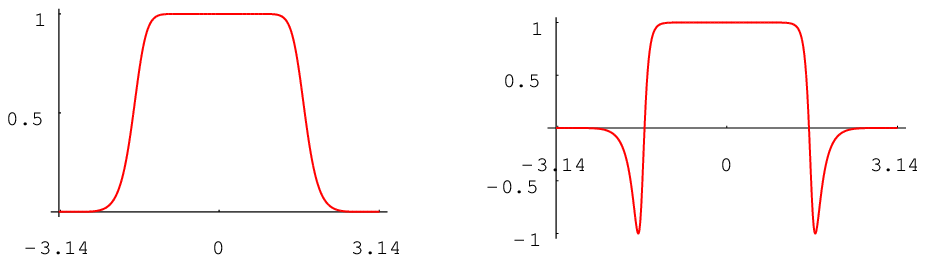}}
\caption{\label{figE} Frequency responses of the maximal
filters $E^0_{12}$ (left) and $E^1_{12}$ (right).}
 For $\delta = 0, 1$
\[E^{\delta}_{12}(e^{i\xi})=\frac{1 + (-1)^{\delta + 1} 
\sqrt{2}\tan^{6}(\xi /2)}
{1 + (-1)^{\delta + 1} \sqrt{2}\tan^{6}(\xi /2)+\tan^{12}(\xi /2)}\] 
 have maximal number of zeros at $\pi$
inside the set of rational linear phase  \qs\ of order $12.$
\end{figure}

\subsection{FIR implementation}
\lb{MF2}
 Following \cite{BEYLKI:1995}, any IIR filter $H(z) =
\frac{P(z)}{Q(z)}$ can be approximated on the unit circle by a FIR 
filter $F(z),$ 
\beq{if1}
F(z) = z^{-N} P(z) \prod_{k = 0}^{n} F_{k}(z^{2^{k}}).
\eeq
The polynomials $F_{k}(z)$ only depend on the zeroes
of $Q(z)$ and on the accuracy sought in the approximation.

$F(z)$ is not an exact QMF anymore, but if for $w$ on the unit circle
\[
\frac{|H(w) - F(w) |}{|H(w)|} \le \epsilon, 
\]
 then, for all $z$
\[
F(z)F(z^{-1}) + F(-z) F(- z^{-1}) = 1 + R(z)
\]
where the FIR filter $R(z)$ satisfies $|R(w)| < 3 \epsilon $ for 
$w$ on the unit circle.  We assumed $\epsilon < \inv{2 \deg Q}.$

As an example, we approximate $E_{12}^{0},$the positive maximally
flat  0-SYM filter
with six zeroes at $-1.$ (See Figure \ref{figE}.) To avoid negative
powers of $z$ in \rf{if1}, we shift
the phase of the filter  by the positive integer $N$
and approximate $z^{N} E^0_{12}(z).$

\begin{figure}[ht]
\epsfxsize = 450pt
\centerline{\epsffile{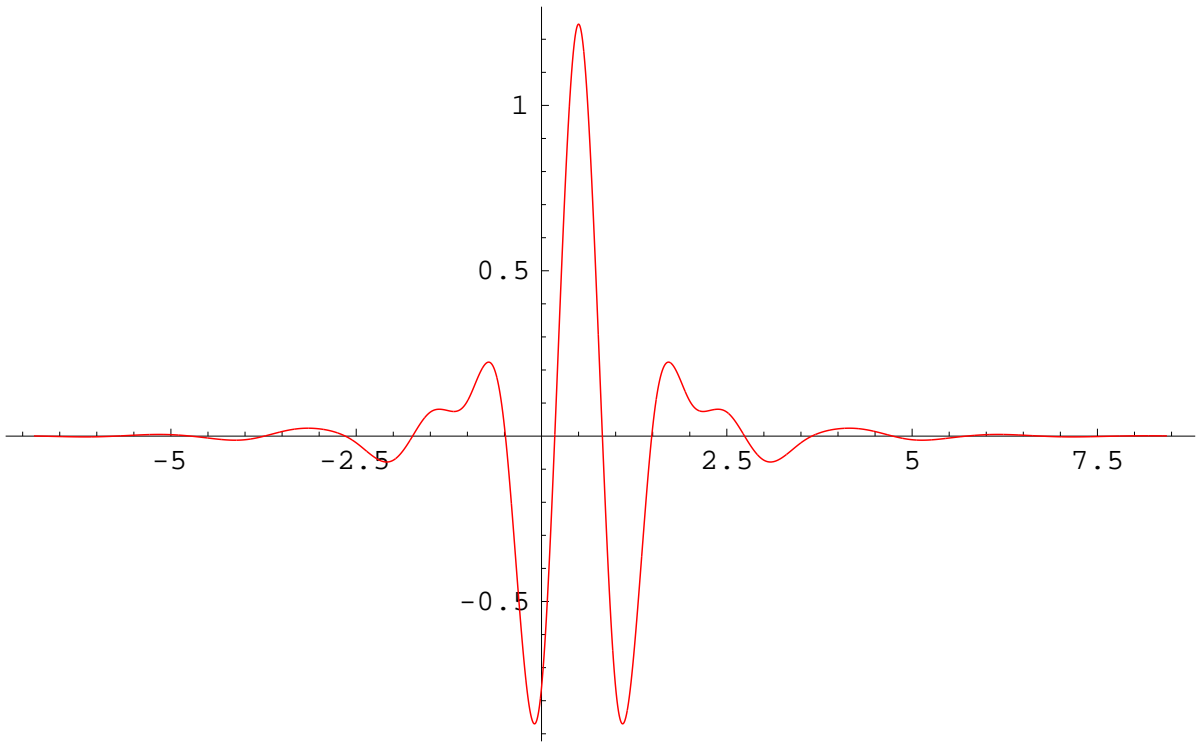}}
\caption{\label{fig11} Wavelet function for FIR approximation
of the maximally flat filter $E^{0}_{12}.$}
\end{figure}

\begin{figure}[htb]
\epsfxsize = 450pt
\centerline{\epsffile{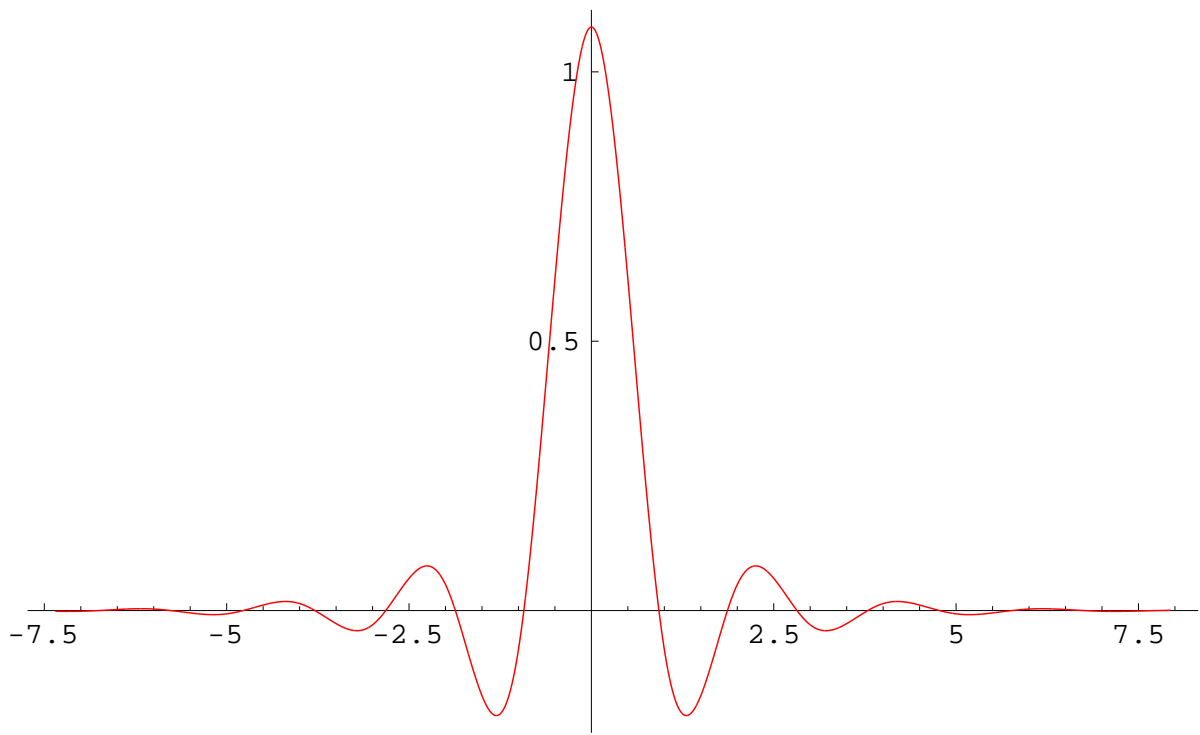}}
\caption{Scaling function for FIR approximation
of the maximally flat filter $E^{0}_{12}.$}
\label{fig12}
\end{figure}

The coefficients of the factors of the approximating filter for
$\epsilon = 10^{-8}$ are listed in Table \ref{tab3}. The factor
$F_0$ equals $1$ because the denominator
of $E^0_{12}$ is an even function. The wavelet and
scaling functions
generated with the approximating filter are plotted in Figures
\ref{fig11} and \ref{fig12}.


\begin{table}[ht] \caption{\lb{tab3} FIR approximation
of a shifted maximal filter with six zeroes at $-1.$The coefficients
listed correspond to the factors in Eq.~\rf{if1}. In this example
$N=0, n = 5,$ and $F_0 = 1.$}
\begin{center}
\begin{tabular}{|c |l | c | l|} \hline
Factor  &  Coefficients  &  Factor & Coefficients\\ [0.5ex]
 \hline\hline
$P$ & -0.0001011263580012439 & $F_2$ & 1.348299677989997e-7 \\
{} & 0.0029296875 & {} & 0.007308809891655256 \\
{} & 0.01818488314800746 & {} & 0.162081739736554 \\
{} & 0.0537109375 & {} & 0.6612186310836452 \\
{} & 0.1156706046299813 & {} & 0.162081739736554 \\
{} & 0.193359375 & {} & 0.007308809891655256 \\
{} & 0.2324912771600248 & {} & 1.348299677989997e-7 \\
{} & 0.193359375 & {} & {} \\
{} & 0.1156706046299813 & $F_3$ & 0.0001276629992294306 \\
{} & 0.0537109375 & {} & 0.03971627520745388 \\
{} & 0.01818488314800746 & {} & 0.920312123586633 \\
{} & 0.0029296875 & {} & 0.03971627520745388 \\
{} & -0.0001011263580012439 & {} & 0.0001276629992294306 \\
{} & {} & {} & {} \\
$F_1$ & -0.00268082617584078 & $F_4$ & 1.925310635034675e-8 \\
{} & 0.6429247852752233 & {} & 0.001585818961857287 \\
{} & -4.433610674839401 & {} & 0.996828323570073 \\
{} & 8.58673343148004 & {} & 0.001585818961857287 \\
{} & -4.433610674839401 & {} & 1.925310635034675e-8 \\
{} & 0.6429247852752233 & {} & {} \\
{} & -0.00268082617584078 & $F_5$ & 2.492072168636633e-6 \\
{} & {} & {} & 0.999995015855663 \\
{} & {} & {} & 2.492072168636633e-6 \\
[0.5ex]
\hline
\end{tabular}
\end{center}
\end{table}

\clearpage
\Section{Conclusion}
\label{conclusion}
The results described in this paper allows one to explicitly construct
all possible linear phase IIR quadrature mirror filters, when the
slope of the phase is zero. Each filter is described using a minimal
set of parameters which are directly related with the zeroes of the
filter. This property permits to construct filters with any preassigned
set of zeroes in the stopband. In particular, we construct, for each
fixed order, the maximally flat filters which generate wavelets with
maximal number of vanishing moments. We also give sufficient
conditions on the parameters to assure that the frequency responses of
the filters are always positive or that the filters have all their
poles on the imaginary axis.  We indicate how to implement these IIR
filters by approximating them by cascades of FIR filters. The FIR
approximation does satisfy the quadrature mirror condition with any
desired accuracy. The factors in the approximation directly depend on
the location of the poles of the IIR filter, what allow one to
introduce the accuracy as part of the filter design process.

\vskip 12pt
\noindent
{\Large \bf Acknowledgments}
\vskip 10pt
\noindent
Some of the results presented can be found in  my PhD thesis. 
I would like to thank my thesis advisor Dr.~R.~Coifman for suggesting the
topic and guiding me through the research.



\bibliographystyle{plain}

%
%
\end{document}